\documentclass[12pt]{amsart}
\usepackage[]{amsfonts}

\def\underset#1#2{\mathrel{\mathop{\kern0pt #2}\limits_{#1}}}

\def\overset#1#2{\mathrel{\mathop{\kern0pt #2}\limits^{#1}}}

\def\couleur(#1 #2 #3)
	{
\special{ps::[end]}}

\def\sqr#1#2{{\vcenter{\vbox{\hrule height.#2pt
             \hbox{\vrule width.#2pt height#1pt \kern#1pt
             \vrule width.#2pt}
             \hrule height.#2pt}}}}

\def\square{\mathchoice\sqr64\sqr64\sqr43\sqr23}			% white square

\def\st{\mathinner{\mkern1mu\raise1pt\hbox{.}				% three points : "such that"
		   \mkern1mu\raise4pt\hbox{.}
		   \mkern1mu\raise1pt\hbox{.}
		 }
         }

\def\bx#1{\setbox1=\hbox{\kern3pt{#1}\kern3pt}				% Make a box. Close it by "}"
 \dimen1=\ht1 \advance\dimen1 by 3pt \dimen2=\dp1 \advance\dimen2 by 3pt
 \setbox1=\hbox{\vrule height\dimen1 depth\dimen2\box1\vrule}%
 \setbox1=\vbox{\hrule\box1\hrule}%
 \advance\dimen1 by .4pt \ht1=\dimen1
 \advance\dimen2 by .4pt \dp1=\dimen2 \box1\relax}

\def\k#1{\kern#1em}
\def\vci{\vrule  width.02em height1.47ex depth-.0ex}				% le 1 en blackboard
\def\11{{\rm\k{.2}\vci\k{-.37}1}}

\newtheorem{Theorem}{Theorem}[section]
\newtheorem{Definition}[Theorem]{Definition}
\newtheorem{Lemma}[Theorem]{Lemma}

\begin{document}
\title{A Carleson type condition for interpolating sequences\\in the Hardy 
spaces of the ball of ${\mathbb{C}}^{n}$.}
\author{Eric Amar}
\address{Universit{\'e} Bordeaux I, 351 Cours de la Lib{\'e}ration\\33405 
Talence France.}
\email{Eric.Amar@math.u-bordeaux1.fr}
\maketitle
\begin{abstract} {
In this work we prove that if $S$ is a dual bounded sequence of points 
in the unit ball ${\mathbb{B}}$ of ${\mathbb{C}}^{n}$ for the Hardy space 
$H^{p}({\mathbb{B}})$, then $S$ is $H^{s}({\mathbb{B}})$ interpolating 
with the linear extension property for any $s\in [1,p[$.\ \par
}\end{abstract}
\section{Introduction.{\hskip 1.8em}}
\setcounter{equation}{0}Let $\displaystyle {\mathbb{B}}$ be the unit ball 
of ${\mathbb{C}}^{n}$ and\textbf{ }$\sigma $ the Lebesgue's measure on 
$\partial {\mathbb{B}}$. As usual we define the Hardy spaces $H^{p}({\mathbb{B}})$ 
as the closure in $L^{p}(\partial {\mathbb{B}})$ of the holomorphic polynomials 
and $\displaystyle H^{\infty }({\mathbb{B}})$ as the algebra of all bounded 
holomorphic functions in $\displaystyle {\mathbb{B}}$.\ \par
{\hskip 1.8em}If $a\in {\mathbb{B}}$ we note $k_{a}(z)$ its Cauchy kernel 
and $k_{a,p}(z)$ its normalized Cauchy kernel in $\displaystyle H^{p}({\mathbb{B}})$:\ 
\par
{\hskip 3.6em}$\displaystyle k_{a}(z):=\frac{\displaystyle 1}{\displaystyle (1-\overline{a}\cdot z)^{n}},\ 
k_{a,p}(z):=\frac{\displaystyle k_{a}}{\displaystyle \displaystyle \left\Vert{\displaystyle 
k_{a}}\right\Vert _{p}}.$\ \par
We know that $\left\Vert{k_{a}}\right\Vert _{p}=c(a,p)(1-\left\vert{a}\right\vert ^{2})^{-n/p'}$ 
with $\frac{1}{p}+\frac{1}{p'}=1$, $0<\alpha \leq c(a,p)\leq \beta $ and 
$\alpha ,\beta $ independent of $a\in {\mathbb{B}}$ and of $1\leq p\leq \infty $.\ 
\par
{\hskip 1.8em}We shall need some definitions.\ \par
\begin{Definition} We say that the sequence $S\subset {\mathbb{B}}$ is 
dual bounded in $\displaystyle H^{p}({\mathbb{B}})$ if a dual system $\{\rho _{a}\}_{a\in S}\subset H^{p}({\mathbb{B}})$ 
for $\displaystyle \{k_{a,p'}\}_{a\in S}$ exists and if this sequence 
is bounded in $\displaystyle H^{p}({\mathbb{B}})$, i.e.\ \par
{\hskip 1.8em}$\displaystyle \exists C>0\st \forall a\in S,\ \displaystyle \left\Vert{\displaystyle 
\rho _{a}}\right\Vert _{p}\leq C,\ \forall a,b\in S,\ \displaystyle \left\langle{\displaystyle 
\rho _{a},\ k_{b,p'}}\right\rangle =\delta _{a,b}.$\ \par
\end{Definition}
\begin{Definition} We say that a sequence $S\subset {\mathbb{B}}$ is $\displaystyle 
H^{p}({\mathbb{B}})$ interpolating for $1\leq p<\infty $, $S\in IH^{p}({\mathbb{B}})$, 
with interpolating constant $C_{I}>0$ if \ \par
{\hskip 3.6em}$\displaystyle \forall \lambda \in \ell ^{p},\ \exists f\in H^{p}({\mathbb{B}})\st \forall 
a\in S,\ f(a)=\lambda _{a}\displaystyle \left\Vert{\displaystyle k_{a}}\right\Vert 
_{p'}$ and $\left\Vert{f}\right\Vert _{p}\leq C_{I}\left\Vert{\lambda }\right\Vert 
_{p}$.\ \par
We say that $S\subset {\mathbb{B}}$ is $\displaystyle H^{\infty }({\mathbb{B}})$ 
interpolating, $S\in IH^{\infty }({\mathbb{B}})$, with interpolating constant 
$C_{I}$ if\ \par
{\hskip 3.6em}$\displaystyle \forall \lambda \in \ell ^{\infty },\ \exists f\in H^{\infty }(\sigma 
)\st \forall a\in S,\ f(a)=\lambda _{a}$ and $\left\Vert{f}\right\Vert _{\infty }\leq C_{I}\left\Vert{\lambda }\right\Vert 
_{\infty }$.\ \par
\end{Definition}
{\hskip 1.8em}Clearly if $S$ is $\displaystyle H^{p}({\mathbb{B}})$ interpolating, 
then $S$ is dual bounded in $\displaystyle H^{p}({\mathbb{B}})$. In one 
variable, L. Carleson~\cite{CarlInt58} proved the converse for $p=\infty $ 
and H. Shapiro \& A. Shields~\cite{ShapShields61} did the same for $\displaystyle 
H^{p}({\mathbb{D}}),\ 1\leq p<\infty $.\ \par
{\hskip 1.8em}In~\cite{amIntInt06} we proved, in the general setting of 
uniform algebras, that if $S$ is $H^{p}$ interpolating then $S$ is $H^{s}$ 
interpolating for any $s$ such that $1\leq s\leq p$ provided that some 
structural hypotheses are verified. This is valid here, in the case of 
the ball.\ \par
\begin{Definition} We say that the $\displaystyle H^{p}({\mathbb{B}})$ 
interpolating sequence $S$ has the linear extension property (LEP) if 
there is a bounded linear operator $E\ :\ \ell ^{p}{\longrightarrow}H^{p}({\mathbb{B}})$ 
such that $\left\Vert{E}\right\Vert <\infty $ and    $\displaystyle \forall \lambda \in \ell ^{p},\ E\lambda $ 
interpolates the sequence $\lambda $ in $\displaystyle H^{p}({\mathbb{B}})$ 
on $S$.\ \par
\end{Definition}
{\hskip 1.8em}In~\cite{amLinExt06} we proved by functional analytic methods 
still in the general setting of uniform algebras, that if $S$ is dual 
bounded in $H^{p}$ then for any $s$ such that $1\leq s\leq p$, $S$ is 
$H^{s}$ interpolating with the LEP, \textsl{provided that }$p=\infty $\textsl{ 
or }$p\leq 2$. Hence it remains a gap, the values of $p\in ]2,\infty [$.\ 
\par
{\hskip 1.8em}The aim of this work is to fill this gap  in the special 
case of the ball.\ \par
\begin{Theorem} \label{extBall21}Let $S\subset {\mathbb{B}}$ be dual bounded 
 in $\displaystyle H^{p}({\mathbb{B}})$, then $S$ is $\displaystyle H^{s}({\mathbb{B}})$ 
interpolating with the LEP, for any $s<p$.\ \par
\end{Theorem}
\section{Carleson sequences.{\hskip 1.8em}}
\setcounter{equation}{0}Remember that $k_{a,q}$ is the normalized reproducing 
kernel for the point $a\in {\mathbb{B}}$ in $\displaystyle H^{q}({\mathbb{B}})$.\ 
\par
\begin{Definition} We say that the sequence $S\subset {\mathbb{B}}$ is 
a Carleson sequence if, for any $q$ such that $1\leq q<\infty $, we have\ 
\par
{\hskip 3.6em}$\displaystyle \exists D_{q}>0,\ \forall \mu \in \ell ^{q},\ \displaystyle \left\Vert{\displaystyle 
\displaystyle \sum_{a\in S}^{}{\mu _{a}k_{a,q}}}\right\Vert _{q}\leq D_{q}\displaystyle 
\left\Vert{\displaystyle \mu }\right\Vert _{q}.$\ \par
\end{Definition}
{\hskip 1.8em}In fact, up to the duality $\displaystyle H^{p}({\mathbb{B}})-H^{q}({\mathbb{B}})$, 
it is proved by H{\"o}rmander~\cite{HormPSH67} that this condition, for 
a $q>1$, is equivalent to the fact that the measure $\displaystyle \chi :=\displaystyle \sum_{a\in S}^{}{(1-\displaystyle \left\vert{\displaystyle 
a}\right\vert ^{2})^{n}\delta _{a}}$ is a Carleson measure, hence if the 
condition is true for a $q>1$ it is true for all $q$'s.\ \par
{\hskip 1.8em}P. Thomas~\cite{Thomas87} proved that if the sequence $S$ 
is $\displaystyle H^{p}({\mathbb{B}})$ interpolating for a $p\geq 1$, 
then $\gamma $ is a Carleson measure, hence $S$ is a Carleson sequence. 
As a corollary of this result we have\ \par
\begin{Lemma} \label{extBall23}If $S$ is dual bounded in $\displaystyle 
H^{p}({\mathbb{B}})$ for a $p\geq 1$ then $S$ is a Carleson sequence.\ 
\par
\end{Lemma}
{\hskip 1.8em}proof:\ \par
the dual system $\{\rho _{a}\}_{a\in S}\subset H^{p}({\mathbb{B}})$ exists 
and is bounded in $\displaystyle H^{p}({\mathbb{B}})$. Let\ \par
{\hskip 1.8em}$\displaystyle \rho _{a,1}:=\rho _{a}k_{a,p'}$ with $p'$ 
the conjugate exponent for $p$. Then we have\ \par
{\hskip 3.6em}$\displaystyle \forall a\in S,\ \displaystyle \left\Vert{\displaystyle \rho _{a,1}}\right\Vert 
_{1}\leq \displaystyle \left\Vert{\displaystyle \rho _{a}}\right\Vert 
_{p}\displaystyle \left\Vert{\displaystyle k_{a,p'}}\right\Vert _{p'}\leq 
\displaystyle \left\Vert{\displaystyle \rho _{a}}\right\Vert _{p}\leq 
C.$\ \par
Moreover $\displaystyle \rho _{a,1}(b)=\rho _{a}(b)k_{a,p'}(b)=\delta _{ab}\rho _{a}(a)k_{a,p'}(a)=\delta 
_{ab}(1-\displaystyle \left\vert{\displaystyle a}\right\vert ^{2})^{-n}$. 
Hence $\{\rho _{a,1}\}_{a\in S}$ is a dual system for $\{k_{a,\infty }\}_{a\in S}$ 
and this means that $S$ is also dual bounded in $\displaystyle H^{1}({\mathbb{B}})$. 
But then it is clear that a dual bounded sequence in $\displaystyle H^{1}({\mathbb{B}})$ 
is $\displaystyle H^{1}({\mathbb{B}})$ interpolating and we can apply 
Thomas' theorem to conclude. \hfill$\square$\ \par
\section{The main result.{\hskip 1.8em}}
\setcounter{equation}{0}We are now in position to prove the theorem~\ref{extBall21}.\ 
\par
We shall need the following lemma which was proved for the Poisson Szeg{\"o} 
kernel in the ball in~\cite{AmarBonami}. We put on $\partial {\mathbb{B}}$ 
the pseudo-distance $\delta (\zeta ,z):=\left\vert{1-\overline{\zeta }\cdot z}\right\vert $ 
and set $D$ its constant in the quasi triangular inequality, i.e.\ \par
{\hskip 1.8em}$\displaystyle \forall a,b\in \partial {\mathbb{B}},\ \delta (a,b)\leq D(\delta (a,c)+\delta 
(c,b))$.\ \par
\begin{Lemma} \label{extBall24}The kernel $\displaystyle K_{t}(\zeta ,z):=\frac{\displaystyle t^{n(p-1)}}{\displaystyle \displaystyle 
\left\vert{\displaystyle 1-(1-t)\overline{\zeta }\cdot z}\right\vert ^{np}}$ 
 verifies for $p>1$:\ \par
(H2)     $\displaystyle \displaystyle \left\vert{\displaystyle K_{t}(\zeta ,z)}\right\vert \leq 
Ct^{\alpha }(\delta (\zeta ,z)+t)^{-\alpha -n}$, with $\alpha =n(p-1)>0$.\ 
\par
(H3)   For any $\zeta ,z_{0},z,t$ such that $\delta (z_{0},z)\leq \frac{1}{2D}(t+\delta (\zeta ,z_{0}))$ 
we have\ \par
{\hskip 4.2em}$\displaystyle \displaystyle \left\vert{\displaystyle K_{t}(\zeta ,z)-K_{t}(\zeta ,z_{0})}\right\vert 
\leq C\delta (z_{0},z)^{1/2}(\delta (\zeta ,z_{0})+t)^{-n-1/2}.$\ \par
\end{Lemma}
{\hskip 1.8em}Proof\ \par
the condition is for small $t$ and we have, with $\alpha =n(p-1)$\ \par
{\hskip 3.9em}$\displaystyle K_{t}(\zeta ,z)\simeq \frac{\displaystyle t^{n(p-1)}}{\displaystyle (t+\delta 
(\zeta ,z))^{np}}=t^{\alpha }((t+\delta (\zeta ,z))^{-n-\alpha },$\ \par
and the condition (H2).\ \par
{\hskip 1.8em}For (H3) we have\ \par

\begin{equation} 
\displaystyle \left\vert{\displaystyle K_{t}(\zeta ,z_{0})-K_{t}(\zeta 
,z)}\right\vert \simeq \frac{\displaystyle t^{n(p-1)}}{\displaystyle (t+\delta 
(\zeta ,z_{0}))^{np}(t+\delta (\zeta ,z))^{np}}\displaystyle \left\vert{\displaystyle 
(t+\delta (\zeta ,z_{0}))^{np}-(t+\delta (\zeta ,z))^{np}}\right\vert 
.\label{extBall22}
\end{equation} \ \par
The function $x{\longrightarrow}f(x):=x^{np}$ verifies\ \par
{\hskip 3.9em}$\left\vert{f(a)-f(b)}\right\vert =\left\vert{b-a}\right\vert \left\vert{f'(c)}\right\vert 
=\left\vert{b-a}\right\vert np\left\vert{c}\right\vert ^{np-1},\ c\in 
]a,b[$ \ \par
by the mean value property, hence here we get\ \par
{\hskip 1.5em}$\displaystyle \displaystyle \left\vert{\displaystyle (t+\delta (\zeta ,z_{0}))^{np}-(t+\delta 
(\zeta ,z))^{np}}\right\vert \leq \displaystyle \left\vert{\displaystyle 
\delta (\zeta ,z_{0})-\delta (\zeta ,z)}\right\vert \max \ \displaystyle 
\left({\displaystyle (t+\delta (\zeta ,z_{0}))^{np-1},\ (t+\delta (\zeta 
,z))^{np-1}}\right) .$\ \par
We have\ \par
{\hskip 2.1em}$t^{n(p-1)}\leq \min \left({(t+\delta (\zeta ,z_{0}))^{n(p-1)},\ (t+\delta 
(\zeta ,z))^{n(p-1)}}\right) ,$\ \par
because $n(p-1)\geq 0$.\ \par
{\hskip 1.8em}Let $A:=\max \left({(t+\delta (\zeta ,z_{0})),\ (t+\delta (\zeta ,z))}\right) 
,\ B:=\min \left({(t+\delta (\zeta ,z_{0})),\ (t+\delta (\zeta ,z))}\right) 
$ then putting the last two inequalities in~(\ref{extBall22}) we get\ 
\par
{\hskip 1.8em}$\displaystyle \displaystyle \left\vert{\displaystyle K_{t}(\zeta ,z_{0})-K_{t}(\zeta 
,z)}\right\vert \leq \displaystyle \left\vert{\displaystyle \delta (\zeta 
,z_{0})-\delta (\zeta ,z)}\right\vert A^{-1}B^{-n}.$\ \par
By the quasi triangular inequality we have\ \par
{\hskip 3.6em}$\displaystyle \delta (\zeta ,z_{0})\leq D\displaystyle \left({\displaystyle \delta (\zeta 
,z)+\delta (z,z_{0})}\right) {\Longrightarrow}t+\delta (\zeta ,z)\geq 
D^{-1}\displaystyle \left({\displaystyle t+\delta (\zeta ,z_{0})}\right) 
-\delta (z_{0},z),$\ \par
because $D\geq 1$. But we have $\delta (z_{0},z)\leq \frac{1}{2D}(t+\delta (\zeta ,z_{0}))$ 
hence, $\displaystyle t+\delta (\zeta ,z)\geq \frac{\displaystyle 1}{\displaystyle 2D}\displaystyle 
\left({\displaystyle t+\delta (\zeta ,z_{0})}\right) $ so we always have 
$\displaystyle A^{-1}B^{-n}\leq C(t+\delta (\zeta ,z_{0}))^{-n-1}$.\ \par
{\hskip 1.8em}Now we shall take advantage that $d(a,b):=\sqrt{\delta (a,b)}$ 
is a distance on $\partial {\mathbb{B}}$. We have\ \par
{\hskip 3.6em}$\displaystyle \displaystyle \left\vert{\displaystyle \delta (\zeta ,z_{0})-\delta (\zeta 
,z)}\right\vert =\displaystyle \left\vert{\displaystyle d^{2}(\zeta ,z_{0})-d^{2}(\zeta 
,z)}\right\vert =\displaystyle \left\vert{\displaystyle d(\zeta ,z_{0})-d(\zeta 
,z)}\right\vert (d(\zeta ,z_{0})+d(\zeta ,z))$.\ \par
Because $d$ is a distance, $\displaystyle \displaystyle \left\vert{\displaystyle d(\zeta ,z_{0})-d(\zeta ,z)}\right\vert 
\leq d(z_{0},z)=\delta ^{1/2}(z_{0},z)$. On the other hand\ \par
{\hskip 1.8em}$\displaystyle d(\zeta ,z)\leq d(\zeta ,z_{0})+d(z_{0},z){\Longrightarrow}d(\zeta ,z_{0})+d(\zeta 
,z)\leq 2d(\zeta ,z_{0})+d(z_{0},z)\leq C\displaystyle \sqrt{\displaystyle 
t+\delta (\zeta ,z_{0})},$\ \par
again with $\delta (z_{0},z)\leq \frac{1}{2D}(t+\delta (\zeta ,z_{0}))$. 
So finally we get\ \par
{\hskip 1.8em}$\displaystyle \displaystyle \left\vert{\displaystyle K_{t}(\zeta ,z_{0})-K_{t}(\zeta 
,z)}\right\vert \leq C\delta ^{1/2}(z_{0},z)(t+\delta (\zeta ,z_{0}))^{-n-1/2},$\ 
\par
and the lemma. \hfill$\square$\ \par
\subsection{Proof of the theorem.{\hskip 1.8em}}
Suppose that the sequence $S$ is dual bounded in $\displaystyle H^{p}({\mathbb{B}})$, 
i.e. there is a sequence $\{\rho _{a}\}_{a\in S}\subset H^{p}({\mathbb{B}})$ 
such that\ \par
{\hskip 3.6em}$\displaystyle \forall a\in S,\ \displaystyle \left\Vert{\displaystyle \rho _{a}}\right\Vert 
_{p}\leq C,\ \displaystyle \left\langle{\displaystyle \rho _{a},\ k_{b,p'}}\right\rangle 
=\delta _{ab}$.\ \par
Now fix $s<p$ and let $\nu \in \ell ^{s}$. We have to show that we can 
interpolate the sequence $\nu $ in $\displaystyle H^{s}({\mathbb{B}})$.\ 
\par
{\hskip 1.8em}Let $q$ be such that $\displaystyle \frac{\displaystyle 1}{\displaystyle s}=\frac{\displaystyle 1}{\displaystyle 
p}+\frac{\displaystyle 1}{\displaystyle q}$ and write $\forall a\in S,\ \nu _{a}=\lambda _{a}\mu _{a}$ 
with $\lambda _{a}:=\frac{\nu _{a}}{\left\vert{\nu _{a}}\right\vert }\left\vert{\nu 
_{a}}\right\vert ^{s/p}$ and $\displaystyle \mu _{a}:=\displaystyle \left\vert{\displaystyle \nu _{a}}\right\vert 
^{s/q};$ then $\lambda \in \ell ^{p},\ \mu \in \ell ^{q}$ and $\displaystyle 
\displaystyle \left\Vert{\displaystyle \nu }\right\Vert _{s}=\displaystyle 
\left\Vert{\displaystyle \lambda }\right\Vert _{p}\displaystyle \left\Vert{\displaystyle 
\mu }\right\Vert _{q}$.\ \par
Let $k_{a,q}$ the normalized reproducing kernel for $a$ in $\displaystyle 
H^{q}({\mathbb{B}})$ and set $\displaystyle \gamma _{a}:=\frac{\displaystyle c(a,q)c(a,s')}{\displaystyle c(a,p')c(a,2)^{2}}$which 
is bounded above and below by strictly positive constants independent 
of $a\in {\mathbb{B}}$. In particular we have $\exists C>0,\ \forall a\in {\mathbb{B}},\ \gamma _{a}\leq C.$\ 
\par
{\hskip 1.8em}Let $\displaystyle h:=\displaystyle \sum_{a\in S}^{}{\nu _{a}\gamma _{a}\rho _{a}k_{a,q}}=\displaystyle 
\sum_{a\in S}^{}{\gamma _{a}\lambda _{a}\rho _{a}\mu _{a}k_{a,q}}$ ; $h$ 
takes the "good" values on $S$ :\ \par
{\hskip 3.6em}$\displaystyle \forall b\in S,\ h(b)=\displaystyle \sum_{a\in S}^{}{\gamma _{a}\nu _{a}\rho 
_{a}(b)k_{a,q}(b)}=\gamma _{b}\nu _{b}\rho _{b}(b)k_{b,q}(b)=\nu _{b}\displaystyle 
\left\Vert{\displaystyle k_{b}}\right\Vert _{s'},$\ \par
because\ \par
{\hskip 2.1em}$\displaystyle \gamma _{b}\rho _{b}(b)k_{b,q}(b)=\gamma _{b}\displaystyle \left\Vert{\displaystyle 
k_{b}}\right\Vert _{p'}\frac{\displaystyle \displaystyle \left\Vert{\displaystyle 
k_{b}}\right\Vert _{2}^{2}}{\displaystyle \displaystyle \left\Vert{\displaystyle 
k_{b}}\right\Vert _{q}}=\gamma _{b}c(b,p')\frac{\displaystyle c(b,2)^{2}}{\displaystyle 
c(b,q)}(1-\displaystyle \left\vert{\displaystyle b}\right\vert ^{2})^{-n/s}=\displaystyle 
\left\Vert{\displaystyle k_{b}}\right\Vert _{s'},$\ \par
by use of $\displaystyle \frac{\displaystyle 1}{\displaystyle s}=\frac{\displaystyle 1}{\displaystyle 
p}+\frac{\displaystyle 1}{\displaystyle q}$.\ \par
{\hskip 1.8em}Moreover $h$ depends linearly on $\nu $.\ \par
{\hskip 1.8em}It remains to estimate its $\displaystyle H^{s}({\mathbb{B}})$ 
norm.\ \par
{\hskip 1.8em}We have, using H{\"o}lder inequalities\ \par
{\hskip 3.6em}$\displaystyle \displaystyle \left\vert{\displaystyle h}\right\vert \leq C\displaystyle 
\left({\displaystyle \displaystyle \sum_{a\in S}^{}{\displaystyle \left\vert{\displaystyle 
\lambda _{a}}\right\vert ^{p}\displaystyle \left\vert{\displaystyle \rho 
_{a}}\right\vert ^{p}}}\right) ^{1/p}\displaystyle \left({\displaystyle 
\displaystyle \sum_{a\in S}^{}{\displaystyle \left\vert{\displaystyle 
\mu _{a}}\right\vert ^{p'}\displaystyle \left\vert{\displaystyle k_{a,q}}\right\vert 
^{p'}}}\right) ^{1/p'}.$\ \par
{\hskip 1.8em}Let $\displaystyle g:=\displaystyle \left({\displaystyle \displaystyle \sum_{a\in S}^{}{\displaystyle 
\left\vert{\displaystyle \lambda _{a}}\right\vert ^{p}\displaystyle \left\vert{\displaystyle 
\rho _{a}}\right\vert ^{p}}}\right) ^{1/p}$, we have\ \par
{\hskip 3.6em}$\displaystyle \displaystyle \left\Vert{\displaystyle g}\right\Vert _{p}^{p}=\displaystyle 
\int_{\partial {\mathbb{B}}}^{}{\displaystyle \sum_{a\in S}^{}{\displaystyle 
\left\vert{\displaystyle \lambda _{a}}\right\vert ^{p}\displaystyle \left\vert{\displaystyle 
\rho _{a}}\right\vert ^{p}}\,d\sigma }=\displaystyle \sum_{a\in S}^{}{\displaystyle 
\left\vert{\displaystyle \lambda _{a}}\right\vert ^{p}\displaystyle \left\Vert{\displaystyle 
\rho _{a}}\right\Vert _{p}^{p}}\leq C\displaystyle \left\Vert{\displaystyle 
\lambda }\right\Vert _{p}^{p},$\ \par
hence $g\in L^{p}(\sigma )$.\ \par
{\hskip 1.8em}So to have $h$ in $L^{s}(\sigma )$ it suffices to prove 
that\ \par
{\hskip 3.6em}$\displaystyle \displaystyle \left({\displaystyle \displaystyle \sum_{a\in S}^{}{\displaystyle 
\left\vert{\displaystyle \mu _{a}}\right\vert ^{p'}\displaystyle \left\vert{\displaystyle 
k_{a,q}}\right\vert ^{p'}}}\right) ^{1/p'}\in L^{q}({\mathbb{B}}).$\ \par
{\hskip 1.8em}So let $\displaystyle f:=\displaystyle \sum_{a\in S}^{}{\displaystyle \left\vert{\displaystyle 
\mu _{a}}\right\vert ^{p'}\displaystyle \left\vert{\displaystyle k_{a,q}}\right\vert 
^{p'}}$ we shall show that $f\in L^{q/p'}$ ($q>p'{\Longrightarrow}\frac{q}{p'}>1$).\ 
\par
{\hskip 1.8em}From $\displaystyle k_{a,q}=\frac{\displaystyle (1-\displaystyle \left\vert{\displaystyle 
a}\right\vert ^{2})^{n/q'}}{\displaystyle (1-\overline{a}\cdot z)^{n}}$, 
we get\ \par
{\hskip 2.1em}$\displaystyle \displaystyle \left\vert{\displaystyle k_{a,q}}\right\vert ^{p'}=\frac{\displaystyle 
(1-\displaystyle \left\vert{\displaystyle a}\right\vert ^{2})^{\displaystyle 
\frac{\displaystyle np'}{\displaystyle q'}}}{\displaystyle \displaystyle 
\left\vert{\displaystyle 1-\overline{a}\cdot z}\right\vert ^{np'}}=(1-\displaystyle 
\left\vert{\displaystyle a}\right\vert ^{2})^{\displaystyle -\frac{\displaystyle 
np'}{\displaystyle q}}\frac{\displaystyle (1-\displaystyle \left\vert{\displaystyle 
a}\right\vert ^{2})^{n(p'-1)}}{\displaystyle \displaystyle \left\vert{\displaystyle 
1-\overline{a}\cdot z}\right\vert ^{np'}}(1-\displaystyle \left\vert{\displaystyle 
a}\right\vert ^{2})^{n}.$\ \par
{\hskip 1.8em}Now the measure $\displaystyle \chi :=\displaystyle \sum_{a\in S}^{}{(1-\displaystyle \left\vert{\displaystyle 
a}\right\vert ^{2})^{n}\delta _{a}}$ is Carleson by lemma~\ref{extBall23} 
and the kernel\ \par
{\hskip 2.1em}$\displaystyle K(a,z):=\frac{\displaystyle (1-\displaystyle \left\vert{\displaystyle 
a}\right\vert ^{2})^{n(p'-1)}}{\displaystyle \displaystyle \left\vert{\displaystyle 
1-\overline{a}\cdot z}\right\vert ^{np'}}$\ \par
verifies (H2) and (H3) by lemma~\ref{extBall24}, then we can apply the 
results in~\cite{AmarBonami} : if $\chi $ is a Carleson measure and $\alpha \in L^{r}(\chi )$ 
then the balayage of the measure $\alpha \,d\chi $ by the kernel $K(a,z)$ 
is in the space $L^{r}(\sigma )$.\ \par
{\hskip 1.8em}Here, let $\eta _{a}(z)$ be a smooth function with support 
in the hyperbolic ball centered at $a$ and with radius $r>0$ so small 
that these balls are disjoint for $a\in S$ and such that $\eta _{a}(a)=1$. 
Let\ \par
{\hskip 2.1em}$\displaystyle \alpha :=\displaystyle \sum_{a\in S}^{}{\displaystyle \left\vert{\displaystyle 
\mu _{a}}\right\vert ^{p'}(1-\displaystyle \left\vert{\displaystyle a}\right\vert 
^{2})^{\displaystyle -\frac{\displaystyle np'}{\displaystyle q}}\eta _{a}(z)},$\ 
\par
if $\displaystyle \alpha \in L^{q/p'}(\chi )$, we have that its balayage, 
which is precisely $f$, is in $\displaystyle L^{q/p'}(\sigma )$.\ \par
{\hskip 1.8em}But we have\ \par
{\hskip 1.8em}$\displaystyle \displaystyle \int_{}^{}{\displaystyle \left\vert{\displaystyle \alpha 
}\right\vert ^{q/p'}\,d\chi }=\displaystyle \sum_{a\in S}^{}{(1-\displaystyle 
\left\vert{\displaystyle a}\right\vert ^{2})^{n}\displaystyle \left\vert{\displaystyle 
\mu _{a}}\right\vert ^{q}(1-\displaystyle \left\vert{\displaystyle a}\right\vert 
^{2})^{-n}}=\displaystyle \sum_{a\in S}^{}{\displaystyle \left\vert{\displaystyle 
\mu _{a}}\right\vert ^{q}}<\infty ,$\ \par
and we are done. \hfill$\square$\ \par

\bibliographystyle{/usr/share/texmf/bibtex/bst/base/plain}

\end{document}